\documentclass[11pt,a4paper,reqno]{amsart}
\usepackage{amsmath}
\usepackage{amssymb,latexsym}
\usepackage[dvips]{graphicx}
\usepackage{color}
\usepackage{cite}
\usepackage{amsthm}

\newtheorem{theorem}{Theorem} [section]
\newtheorem{lemma}{Lemma} [section]
\newtheorem{proposition}{Proposition} [section]

\newtheorem{definition}{Definition} [section]

\newtheorem{remark}{Remark}[section]

\let\ssection=\section\renewcommand{\section}{\setcounter{equation}{0}\ssection}
\begin{document}
\date{}

\address{M. Darwich: Universit\'e Fran\c{c}ois rabelais de Tours, Laboratoire de Math\'ematiques
et Physique Th\'eorique, UMR-CNRS 7350, Parc de Grandmont, 37200
Tours, France} \email{Mohamad.Darwich@univ-tours.fr}
\title{On the invariant  measures for 
 the Ostrovsky  equation.}
\author{Darwich Mohamad.}
\begin{abstract}
In this paper, we construct
 invariant measures for the Ostrovsky equation associated with conservation laws. On the other hand, 
we prove  the local well- posedness of the initial value problem for the periodic Ostrovsky equation with initial data in $H^{s}(\mathbb{T})$
 for $s>-\frac{1}{2}$. 
 \end{abstract}

\maketitle



\section{Introduction}
 In this paper, we construct an invariant measure for a dynamical system defined by the Ostrovsky 
equation ($\text{Ost}$)

\begin{equation}\label{Ost}
\left\{\begin{array}{l}
\partial_tu - u_{xxx} + \partial_{x}^{-1}u +uu_x=0,\\
u(0,x)= u_0(x).
\end{array}\right.\end{equation}
associated to the conservation of the energie.
The operator $ \partial_x^{-1}$ in the equation denotes a certain antiderivative with respect to the variable $x$ defined for $0$-mean value periodic
function the Fourier transform by 
$\widehat{(\partial_{x}^{-1}f)} = \frac{\hat f (\xi)}{i\xi}$.\\

Invariant measure play an important role in the theory of dynamical systems (DS). 
It is well known that the whole ergodic theory is based on this concept. On the other hand, they are necessary in various physical considerations.\\
Note that, one the well-known applications of invariant measures in the theory of dynamical  
is the Poincar\'{e} recurrence theorem : every flow which
preserves a finite measure has the returning property modulo a set of measure zero.\\
Recently several papers(\cite{ARSEN},\cite{Zhidkov},\cite{Zhidkov1}) have been published on invariant measures for dynamical system generated by nonlinear partial differentiel equations. \\
In \cite{Zhidkov2} an infinite series of invariant measure associated with a higher conservation laws are
constructed for the one-dimensional Korteweg de Vries (KdV) equation:
$$
u_t+uu_x+u_{xxx}=0,
$$ by Zhidkov. In particular, invariant
measure associated to the conservation of the energie are constructed for this equation.\\
 Equation \ref{Ost} is a perturbation of the  Korteweg de Vries (KdV) equation
with a nonlocal term and was deducted by Ostrovskii \cite{Ostrovskii} as a
model for weakly nonlinear long waves, in a rotating frame of reference, to describe the propagation of surface waves in the
ocean. \\

We will construct invariant measures associated to the conservation 
of the Hamiltonian:
$$H(u(t)) = \frac{1}{2} \int (u_x)^2 + \frac{1}{2}\int (\partial_x^{-1} u)^2 - \frac{1}{6}\int u^3.$$

The paper is organized as follows. In Section \ref{Notmesures} the
basic notation is introduced and the basic results are formulated. In Section \ref{Invariantemeasure} the
invariant measure which corresponds to the conservation of the Hamiltonian is constructed.\\
 In Section \ref{wellposednessinBourgain} we will prove the local  well-posedness for our equation in $H^{s}$, $s > -\frac{1}{2}$.

\section{Notations and main results}\label{Notmesures}
We will use $C$ to denote various time independent constants, usually depending only upon $s$.
In case a constant depends upon other quantities, we will try to make it explicit.
 We use $A \lesssim B$ to denote an estimate of the form $A\leq C B$. similarly,
 we will write $A\sim B$ to mean $A \lesssim B$  and $B \lesssim A$. We  writre
 $\langle \cdot\rangle :=(1+|\cdot|^2)^{1/2}\sim 1+|\cdot| $. The notation $a^+$ denotes
$a+\epsilon$ for an arbitrarily small $\epsilon$. Similarly $a-$ denotes $a-\epsilon$.
Let$$ L^2_0 = \{ u \in L^2; \int_{\mathbb{T}} u dx =0\}.$$
On the circle, the Fourier transform is
defined as
$$ \hat{f}(n)=\frac{1}{2\pi}\int_{\mathbb{T}}f(x)\exp(-inx)dx.$$
We  introduce the zero mean-value Sobolev spaces   $H^{s}$  defined by~:
\begin{equation}\label{Esp1}
H_0^{s}=:\{ u \in \mathcal S^{'}(\mathbb{T}); ||u||_{H^{s}}
<+\infty~\text{and}~ \int_{_{\mathbb{T}}} u dx=0\},
\end{equation}
where,
\begin{equation}\label{Ns}
||u|| _{H_0^{s}}= (2\pi)^{\frac{1}{2}}||{\langle.\rangle}^{s}\hat u||_{l^2_n},
\end{equation}
and $X^{s,\frac{1}{2}}$ by 
$$
\{ u \in \mathcal S^{'}(\mathbb{T}); ||u||_{X^{s,\frac{1}{2}}}:=||{\langle n\rangle}^{s}\langle \tau + n^3 - \frac{1}{n}\rangle \hat u||_{l^2_nL^2_\tau}< \infty\}.
$$

Let $$Y^{s} =:\{ u \in \mathcal S^{'}(\mathbb{T}); ||u||_{Y^{s}}
<+\infty \},$$
where
$$
||u|| _{Y^{s}}= ||u|| _{X^{s,\frac{1}{2}}} +||{\langle n\rangle}^{s}\hat u(n,\tau||_{l^2_nL^1_\tau}.
$$

We will briefly remind the general construction of a Gaussian measure on a Hilbert space. Let $X$ be a Hilbert space, and $\{e_k\} $ be the orthonormal basis in $X$ which consists of eigenvectors of some operator $S=S^* > 0$ with corresponding  eigenvalues $0 < \lambda_1 \leq \lambda_2 \leq \lambda_3 ....\leq \lambda_k \leq ...$ We call a set $M \subset X$ a cylindrical
set iff:
$$
M =\{x \in X; [(x,e_1),(x,e_2),...(x,e_r)]\in F\}
$$
for some Borel $F \subset \mathbb{R}^{r}$, and some integer $r$. We define the measure $w$ as follows:
\begin{equation}\label{Gaussian}
w(M) = (2\pi)^{-\frac{r}{2}}\prod_{j=1}^{r}\lambda_j^{\frac{1}{2}}\int_{F}e^{-\frac{1}{2}\sum_{j=1}^{r}\lambda_{j}y_j^2}dy.
\end{equation}
One can easily verify that the class $\mathbb{A}$ of all cylindrical sets is an algebra on which the function $w$ is additive.
The function $w$ is called the centered Gaussian measure  on $X$ with the correlation operator $S^{-1}$.
\begin{definition}
The measure $w$ is called a  countably additive measure on an algebra $\mathbb{A}$ if  $\lim_{n \rightarrow +\infty}(A_n)=0$ for any $A_n \in \mathbb{A}$($n=1,2,3...$) for which $A_1 \supset A_2 \supset A_3\supset.....\supset A_n\supset...$ and $\bigcap_{n=1}^{\infty}A_n = \phi$
\end{definition}
Now we give the following Lemma:
\begin{lemma}\label{lemmacountablyadditive}
The measure $w$ is countably additive on the algebra $\mathbb{A}$ iff $S^{-1}$ is an operator of trace class, i.e iff 
$\sum_{k=1}^{+\infty}\lambda_k^{-1} < +\infty $.
\end{lemma}
Now we present some definitions related to invariant  measure :
\begin{definition}
Let M be a complete separable metric space and let a function $h: \mathbb{R} \times M\longmapsto M$ for any fixed $t$ be a homeomorphism of the space $M$ into itself satisfying the properties:\\
\begin{enumerate}
\item $h(0,x)=x$ for any $x \in M$.
\item $h(t,h(\tau,x))=h(t+\tau,x)$ for any $t,\tau \in \mathbb{R}$ and $x \in M$.
\end{enumerate}
Then, we call the function $h$ a dynamical system with the space $M$. If $\mu$ is a Borel measure defined on the phase 
space $M$ and $\mu(\Omega)=\mu(h(\Omega,t))$for an arbitrary Borel set $\Omega \subset M$ and for all $t \in \mathbb{R}$, 
then it is called an invariant measure for the dynamical system $h$.
\end{definition}
Let us now state our results:
\begin{theorem}\label{thexistence}

Let $s>-1/2$,  and $\phi \in H_0^{s}$. Then
 there exists a time $T=T(||\phi||_{H_0^{s}})>0$ and a unique solution $u$ of  (\ref{Ost}) 
in $ C([0,T],H_0^{s})\cap Y^s$ and  the map $\phi\longmapsto u$ is $C^\infty$ from $H_0^{s}$ to $C([0,T],H_0^{s})$.
$\hfill{\Box}$
\end{theorem}

\begin{theorem}\label{theoremmeasure}
 Let $\phi \in L_0^2$, then  the Problem \ref{Ost} is global well-posedness in $L^2$ and
  the Borel measure $\mu$  on $L^2$ defined for any Borel set $\Omega \subset L^2$ by the rule
 $$
 \mu(\Omega)=\int_{\Omega}e^{-g(u)}dw(u)
 $$
 where $w$ is the centered Gaussian measure corresponding to the correlation operator $S^{-1}= (-\Delta + \Delta^{-1})^{-1}$, and $g(u)= \frac{1}{3}\int u^3 dx$ 
the nonlinear term of the Hamiltonian 
 is an invariant measure for (\ref{Ost}).
\end{theorem}
 

\section{Invariance of Gibbs measure}\label{Invariantemeasure}
In this section, we construct an invariant measure to Equation \ref{Ost}  with respect to the conservation of the Hamiltonian. Let us first 
  present result on invariant measures for systems of autonomous ordinary differential equations. Consider the following system of ordinary differential
 equations:
\begin{equation}\label{ordinary}
\dot{x}=b(x),
\end{equation}
where $x(t):\mathbb{R}\longmapsto \mathbb{R}^{n}$ is an unknown
 vector-function and $b(x):\mathbb{R}^{n}\longmapsto \mathbb{R}^{n}$ is a continuously differentiable map. Let $h(t,x)$ 
be the corresponding function (`` dynamical system'') from $\mathbb{R}\times\mathbb{R}^{n}$ into $\mathbb{R}^{n}$ transforming any $t \in \mathbb{R}$ 
and $x_0 \in \mathbb{R}^{n}$ into the solution $x(t)$, taken at the moment of time $t$, of the above system supplied with the initial data $x(0)=x_0$. 
\begin{theorem}\label{equivalence}
 Let $P(x)$ be a continuously differentiable function from $\mathbb{R}^{n}$ into $\mathbb{R}$. For the Borel measure 
$$
\nu(\Omega)=\int_{\Omega}P(x)dx
$$
to be invariant for the function $h(t,x)$ in the sense that $\nu(h(t,\Omega))=\nu(\Omega)$ for any bounded domain $\Omega$ and for any $t$, it is sufficient and 
necessary that
$$
\displaystyle{\sum_{i=1}^{n}}\frac{\partial}{\partial x_i}(P(x)b_i(x)) = 0,
$$
for all $x \in \mathbb{R}^{n}$.
\end{theorem}



We shall construct an invariant measure for (\ref{Ost}). Let $A >0$, the space $L^2(0,A)$ be real equipped with the scalar product:
$$
(u,v)_{L^2(0,A)} = \int_0^A u\overline{v}dx.
$$
and $J = \frac{\partial}{\partial x}Q$ where the operator $Q$ maps $v^* \in L^2$ into $v \in L^2$ such that
$v^*(g) = (v,g)_{L^2(0,A)}$.
 Finally, let $S = -  \Delta + \Delta^{-1}$. We set $H(u) = \frac{1}{2}(\int (u_x)^2 - 
\int(\partial_x^{-1} u )^2 ) + \frac{1}{3} \int u^3= \frac{1}{2}(Su,u) + g(u)$.
Note that System \ref{Ost} takes the form:
\begin{equation}\label{abstrat}
\left\{\begin{array}{l}
\frac{\partial u}{\partial t}(t) = J \frac{\delta}{\delta u } H(u(t)), t \in \mathbb{R}\\
u(t_0)= \phi \in H^{s},
\end{array}\right.\end{equation}
Let $ e_{2k-1}(x)=\frac{\sqrt 2}{\sqrt A} \text{sin}( \frac{2 \pi n x}{A})$, $ e_{2k} = \frac{\sqrt 2}{\sqrt A} \text{cos} (\frac{2 \pi n x}{A})$
where $ k = 1,2,3...$ Then $(e_k)_{k=1,2,..}$ is an orthonormal basis of the space $L_0^2(0,A)$ consisting of eigenfunctions of the operator $\Delta$
 with the corresponding eigenvalues $0  < \lambda_1 = \lambda_2 < ...< \lambda_{2k-1} = \lambda _{2k} < ...$ Let $P_m$ be  the 
orthogonal projector in $L^2_0$ onto the subspace $L_m = \text{span}\{e_1,...,e_{2m}\}$ and $P_m^{\bot}$ be the orthogonal projector in $L_0^2(0,A)$ onto the 
orthogonal complement $L_m^{\bot}$ to the subspace $L_m$. Let also $ v_i = -\lambda_i + \lambda_i^{-1}$, then $v_i$ are eigenvalues of $S$.\\
Consider  the following problem:

\begin{equation}\label{Ostm}
\left\{\begin{array}{l}
\partial_tu^m - u^m_{xxx} + \partial_{x}^{-1}u^m + P_m(u^mu^m_x)=0,\\
u^m(0,x)=P_m u_0(x).
\end{array}\right.\end{equation}
The existence of $u$ is global in $L^2$ in time ( see later) and
the solution of (\ref{Ostm}) converges to $u$ in $C([0,T], L^2)$ for any fixed $T$, more precisely we have the following lemma:
\begin{lemma}\label{convsolution}
\begin{enumerate}
 \item
The solution $u_m$ of (\ref{Ostm}) converges in $C([0,T], L^2)$ to the solution $u$ of (\ref{Ost}).\item 
For any $\epsilon > 0$, and $T>0$ there exists $ \delta > 0$ such that 
$$
\displaystyle{\text{Max}_{ t \in [t_0-T,t_0+T]}}\lVert u_m(.,t)-v_m(.,t)\rVert_{L^2} < \epsilon,
$$

for any two solutions $u_m$ and $v_m$ of the problem (\ref{Ostm}), satisfying the condition 
$$
\lVert u_m(.,t_0)-v_m(.,t_0)\lVert_{L^2} < \delta.
$$
 \end{enumerate}
\end{lemma}
\proof : 
By the Duhamel formula, $ u -u^m$ satisfies
$$
u(t)-u^m(t)= e^{-itS}(u_0-P_mu_0) - \frac{1}{2} \int_0^t e^{-i(t-t')S}(\partial_x(u^2(t'))-P_m(\partial_x((u^m)^2(t'))))dt'.
$$
We can whrite that $R(t):= \partial_x(u^2(t'))-P_m(\partial_x((u^m)^2(t'))) = \partial_x(u^2-(P_\frac{m}{2}u)^2) + 
P_m\partial_x\bigg((P_\frac{m}{2}u)^2 - u^2\bigg) + P_m\partial_x(u^2 - (u_m)^2).$ Now, using the linear and bilinear estimates 
proved in section \ref{wellposednessinBourgain}, we obtain that
\begin{equation}
 \|u-u_m\|_{Y^s} \lesssim \|u_0-P_mu_0\|_{H^s} +T^{\gamma}\|u-u_m\|_{Y^s}\|u+u_m\|_{Y^s} + \|u-P_{\frac{m}{2}}u\|_{Y^s},
\end{equation}

then $u_m \longrightarrow u $ in $Y^s$, but $Y^s \hookrightarrow L^{\infty}_t{L^2_x}$, this gives the uniform convergence in $L^2$. \\
The proof of part (2) is similar to part (1).\\

By $h_m(u_0,t)$ we denote the function  mapping any $u_0 \in L^2$ and $ t \in \mathbb{R}$ into $u_m(.,t+t_0)$ where $u_m(.,t)$ 
is the solution of the problem
(\ref{Ostm}). It is clear that the function $h_m$ is a dynamical system with the phase space $X^m = \text{span} \{e_1,...e_m\}$.
 In addition, the direct verification shows that $\frac{d}{dt}||u_m(.,t)||_{L^2}^2=0$ and $\int u_m dx =0$.
. For each $m=1,2...$ 
let us consider in the space $X^m$ the centered Gaussian measure $w_m$ with the correlation operator $S^{-1}$. Since $S=S^*$ in $X^m$, the measure $w_m$
is well-defined in $X^m$.
Also, since $g(u) = \frac{1}{3}\int u^3$ is a continuous functional in $X^m$, the following Borel measures
$$
\mu_m(\Omega) = \int_{\Omega}e^{-g(u)}dw_m(u).
$$
(where $\Omega$ is an arbitrary Borel set in $L^2$) are well defined.

\begin{definition}A set $\Pi$ of measures defined on the Borel sets of a topologogical space is called tight if, for each $ \epsilon >0$, there exist a compact set
$K$ such that 
$$
\mu(K) > 1-\epsilon
$$
For all $\mu \in \Pi$.

\end{definition}
 We will use the following theorem:
\begin{theorem}(Prokhorov)
A tight set, $\Pi$, of  measures on the Borel sets of a
metric topological space, $X$, is relatively compact in the sense that for each
sequence,$ P_1 , P_2 , . . . $in $\Pi$ there exists a subsequence that converges to a
probability measure $P$ , not necessarily in $\Pi$, in the sense that
$$
\int gdP_{n_j} \longrightarrow \int gdP
$$
for all bounded continuous integrands. Conversely, if the metric space is
separable and complete, then each relatively compact set is tight.

\end{theorem}

To prove Theorem \ref{theoremmeasure}, we will prove the following Lemma:
\begin{lemma}
$\mu _m$ is an invariant measure for the dynamical system $h_m$ with the phase space $X^m$.
\end{lemma}
\textbf{Proof}: Let us rewrite the system (\ref{Ostm}) for the coefficients $a_k$, where $u^m(t)=\displaystyle{\sum_{k=1}^{k=2m}a_k(t)e_k}$. Let 
$h(a) = H(\displaystyle {\sum_{k=1}^{k=2m}a_ke_k})$ and  $J$ is a skew-symmetric matrix, $(J_m)_{2k-1,2k} = -\frac{2\pi k}{A} = -(J_m)_{2k,2k-1}$(k=1,2,.. m)
then 
the problem take the form
\begin{equation}\label{Osta}
\left\{\begin{array}{l}
a^{\prime}(t) = J_m\nabla_{a} h(a(t)),\\
a_k(t_0)= (u_0,e_k), k=1,2,...2m
\end{array}\right.
\end{equation}
Using Theorem \ref{equivalence}, we can easily verify that the Borel measure:
$$
\mu_m^{\prime}(A)= (2\pi)^{-\frac{2m+1}{2}}\prod_{j=1}^{2m}v_j^{\frac{1}{2}}\int_{A}e^{-\frac{1}{2}\sum_{j=1}^{2m}v_{j}a_j^2-g(\sum_{j=1}^{2m}a_je_j(x))}da,
$$ 
(with $v_j= -\lambda_j+ \lambda_j ^{-1}$ the eigenvalues of $S$) is invariant for the problem (\ref{Osta}). Also, we introduce the measures 
$$
w_m(A) = (2\pi)^{-\frac{2m+1}{2}}\prod_{j=1}^{2m}v_j^{\frac{1}{2}}\int_{A}e^{-\frac{1}{2}\sum_{j=1}^{2m}v_{j}a_j^2}da.
$$
Let $\Omega_m \subset X^m$ and $\Omega_m = \{u \in L^2, u = \displaystyle \sum_{j=1}^{2m}a_je_j, a \in A\}$ where $A \subset \mathbb{R}^{2m}$ is a Borel set. We set $\mu_m(\Omega_m)= \mu_m^{\prime}(A)$.
Since the measure $\mu_m^{\prime}$ is invariant for (\ref{Osta}), the measure $\mu_m$ is invariant for the problem (\ref{Ostm}).\\
Although  the measure is defined on $X^m$, we can define it on the Borel sigma-algebra of $L^2$ by the rule: $\mu_m(\Omega)=\mu_m(\Omega \cap X^m)$. Since the set $\Omega \cap X^m$ is open as a set in $X^m$ for any open set $\Omega \subset L^2$, this procedure is correct.

\begin{lemma}\label{weakly}
$(w_m)_m$ weakly converges to $w$ in $L^2$.
\end{lemma}
\textbf{Proof}:
$S^{-1}$ is an operator of trace since the trace $Tr(S^{-1})= \displaystyle{\sum_{k}v_k^{-1}}=\displaystyle{\sum_{k}}\frac{1}{\frac{1}{\frac{4\pi^2k^2}{A^2}}+\frac{4\pi^2k^2}{A^2}}<+\infty$. 
Thus we can find  a continuous positive function $d(x)$ defined on $(0,\infty)$ 
with the property $\displaystyle{\lim_{x \rightarrow +\infty}} d(x) = +\infty$ such that $\displaystyle{\sum_{k}v_k^{-1}d(\lambda_k)} < +\infty$. 
We define the operator $T = d(S)$, the operator defined by $T(e_k)= d(v_k)e_k$ and let $B = S^{-1}T$. According to the definition of $d(x)$, $Tr (B) < + \infty$.
 Let $R > 0$ and $B_R = \{ u \in L^2, T^{\frac{1}{2}} u \in L^2 \text{and} ||T^{\frac{1}{2}}u|| \leq R\}$, it is clear that the closure of $B_R$ is compact for any $R > 0$.
  Combined the following inequality ( see 
 \cite{Fomin} for the proof)
$$
w_n(\overline{B_R}^{C})= w_n(\{u; (Tu,u)_{L^2} > R\}) \leq \frac{Tr (B)}{ R^2}.
$$ 
with the Prokhorov theorem, this ensure that $(w_n)$ is weakly compact on $L^2$.\\
In view of the definition $w_n(M) \rightarrow w(M)$ for any cylindrical set $M \subset L^2$.(because $w_n(M)=w(M)$ for all sufficiently large $n$). 
Hence, since the extension of a measure from an algebra to a minimal sigma-algebra is unique, we have proved that the sequence $w_n$ 
converges to $w$ weakly in $L^2$
and Lemma \ref{weakly} is proved.
\begin{lemma}\label{infsup}
$\liminf_{m}\mu_m(\Omega) \geq \mu (\Omega) $ for any open set $\Omega \subset L^2$.\\
 $\limsup_{m}\mu_m(K) \leq \mu (K) $ for any closed bounded set $K \subset L^2$.\\
\end{lemma}
\textbf{Proof}: 
 Let $\Omega \subset L^2$ be open and let 
 $B_R = \{ u \in L^2,  ||u||_L^2 < R\}$ for some $R > 0$.\\
Consider 
$\phi(u) : 0 <\phi(u)< 1$  with the support belonging to  $\Omega_R=\Omega \cap B_R$ such that
$$
\int_X \phi(u)e^{-g(u)}dw(u) \geqslant \mu(\Omega_R) - \epsilon.
$$
Then,
\begin{align}
\liminf_{m}\mu_m(\Omega_R) & = \liminf_{m}\int_{\Omega_R}e^{-g(u)}dw_m(u) \geq \liminf_{m}\int \phi(u)e^{-g(u)}dw_m(u)  \nonumber\\
&= \int \phi(u)e^{-g(u)}dw(u) \geq \mu( \Omega_R) -\epsilon.\nonumber
\end{align}
Therefore, due to the arbitrariness of $\epsilon >0$ one has:
$$
\liminf_{m}\mu_m(\Omega) \geq \limsup_{m}\mu_m(\Omega_R) \geq \mu(\Omega_R). 
$$
Taking $R \longrightarrow +\infty$ in this inequality, we obtain the first statement the lemma.\\

Let $K$ be a closed bounded set. Fix $\epsilon >0$. We take a continuous function $\phi \in [0,1]$ such that $\phi(u)=1$ for any $ u\in K$, 
$\phi(u) = 0$ if $\text{dist}(u,K) > \epsilon$ and $\int \phi(u)e^{-g(u)}w(du) < \mu(K) + \epsilon$. Then
\begin{align}
 \limsup_{m}\mu_m(K) &\leq \limsup_{m}\int \phi(u)e^{-g(u)}dw_m(u)  \nonumber\\
&= \int \phi(u)e^{-g(u)}dw(u) \leq \mu(K) +\epsilon,\nonumber
\end{align}
and due to the arbitrariness of $\epsilon >0$, Lemma \ref{infsup} is proved.




\begin{lemma}\label{invariantelemma}
Let $\Omega \subset L^2$ an open  set and $t \in \mathbb{R}$. Then $\mu(\Omega)=\mu(h(\Omega,t))$.
\end{lemma}
\textbf{Proof}: Let $\Omega_1=h(\Omega,t)$. Fix an arbitrary $ t\in \mathbb{R}$, then $\Omega_1$ is open too. 
First, let us suppose that $\mu(\Omega) < \infty$.\\
 Fix an arbitrary $\epsilon >0$,
by Prokhorov Theorem there exists a compact set $K \subset \Omega$ such that $\mu(\Omega \backslash K) < \epsilon$, note that $K_1= h(K,t)$ is a compact set, too, and
 $K_1 \subset \Omega_1$. \\
For any $A \subset L^2$, let $\partial A$ be the boundary of the set $A$ and let
$$
\beta=\min\{dist(K,\partial \Omega);dist(K_1,\partial \Omega_1)\}
$$
(where  $\text{dist}(A,B)= \inf_{x\in A,y\in B}\| x-y\|_{L^2}$). Then, $\beta > 0$. According to Lemma \ref{convsolution}, for any $ z \in K$, there exists $\delta >0$
 such that for any $x,y \in B_{\delta}(z)$ one has $\| h_n(x,t) - h_n(y,t)\|_{L^2} < \frac{\beta}{3}$. 
Lets $\Omega^{\alpha} = \{ q \in \Omega_1; \text{dist}(q,\partial \Omega_1) \geq \alpha \}$ and
$B_{\delta_1}(z_1),...B_{\delta_l}(z_l)$ be a finite covering of the compact set $K$ by these balls and let $B= \bigcup_{i=1}^{l}B_{\delta_i}(z_i)$.\\
Since $ h_n(z_i,t)\longrightarrow h(z_i,t) (n \longrightarrow +\infty)$ for any $i$ we obtain that 
$\text{dist}(h_n(z,t), K_1) < \frac{\beta}{3}$, $\forall z \in B$ and  large $n$.
Thus, $h_n(B,t)$ belongs to a closed bounded subset of $\Omega^{\frac{\beta}{2}}$ for all sufficiently large $n$.\\
 Further, we get by the invariance of $\mu_n$ and Lemma \ref{infsup}
$$
\mu(\Omega) \leq \mu(B) +\epsilon \leq \liminf \mu_n(B) + \epsilon \leq \liminf \mu_n(h_n(B,t)) + \epsilon \leq \mu(\Omega_1) + \epsilon
$$ 
$\bigg(\text{because}~~ \mu_{n}(B) = \mu_n(B \cap X_n) = \mu_n(h_n(B \cap X_n,t)), \text{and}~~ h_n(B \cap X_n,t) \subset h_n(B,t)\bigg)$.
Hence, due to the arbitariness of $\epsilon >0$, we have $ \mu(\Omega) \leq \mu(\Omega_1)$. By analogy $ \mu(\Omega) \geqslant \mu(\Omega_1)$. Thus 
$ \mu(\Omega) = \mu(\Omega_1)$.\\
Now if $\Omega$ is open and $\mu(\Omega) = + \infty$, then we take the sequence 
$$
\Omega^{k} = \Omega \cap \{u \in L^2; \|u\|_{L^2} + \|h(u,t)\| < k \}
$$ 
and set $\Omega_1^k= h(\Omega^k,t)$. Then $\Omega = \cup \Omega^k$ and $ \mu(\Omega^k)= \mu(\Omega^k_1) < \infty$. Taking $ k\longrightarrow +\infty$, we 
obtain the statement of the lemma. 



\section{Well-posedness in $X^{\lowercase{s},\frac{1}{2}}$}\label{wellposednessinBourgain}
In this section, we prove a global wellposedness result for the Ostrovsky equation by
following the idea of Kenig, Ponce, and Vega in \cite{Kenig1}.\\

Our work space is $Y^s$, the completion of functions that are Schwarz in time and
$C^\infty$ in space with norm:
$$
||u|| _{Y^{s}}= ||u|| _{X^{s,\frac{1}{2}}} +||{\langle n\rangle}^{s}\hat u(n,\tau)||_{l^2_nL^1_\tau}
$$
$Y ^s$ is a
slight modification of $X^{s,\frac{1}{2}}$ such that $||u|| _{L_t^{\infty}H_x^{s}}\lesssim||u|| _{Y^{s}}$.\\
 We see that the nonlinear part of the Ostrovsky
equation is $u\partial_{x}u$, and by Fourier transform we write it in frequency as
$$
n\displaystyle{\sum_{n_1 \in \mathbb{\dot{Z}}}\int_{\tau_1 \in \mathbb{R}}\hat{u}(n_1,\tau_1)\hat{u}(n-n_1,\tau-\tau_1)d\tau_1}.
$$
The resonance function is given by:
$$
R(n,n_1)= \tau+m(n) - (\tau_1+m(n_1)-(\tau-\tau_1+m(n-n_1)= 3nn_1(n-n_1)-\frac{1}{n}\big(1-\frac{n^3}{nn_1(n-n_1)}\big)
$$
where $m(n)=n^3-\frac{1}{n}$.\\
Now we have the following lower bound on the resonance function:
\begin{lemma}\label{lemma1}
  If $|n||n_1||n-n_1| \neq 0$, and $\frac{1}{\mid n \mid} < 1$, then:
\begin{equation}\label{R}
|R(n,n_1)| \gtrsim |n||n_1||n-n_1|,
\end{equation}
and
\begin{equation}\label{n2}
|n|^2\leq 2 |nn_1(n-n_1)|.
\end{equation}
\end{lemma}
\textbf{Proof}: (\ref{n2}) is obvious.\\
Now 
 \begin{eqnarray*}
R^{2}(n,n_1)&=&9n^2n_1^2(n-n_1)^2  - 6n_1(n-n_1) + 6n^2 + \frac{1}{n^2}\big(1-\frac{n^3}{n(n_1(n-n_1))}\big)^2 \nonumber\\
&=&n^2n_1^2(n-n_1)^2 + 8n^2n_1^2(n-n_1)^2  - 6n_1(n-n_1) + 6n^2 + \frac{1}{n^2}\big(1-\frac{n^3}{n(n_1(n-n_1))}\big)^2\nonumber\\
&\ge& n^2n_1^2(n-n_1)^2 + 8n^2n_1^2(n-n_1)^2  - 6n_1(n-n_1)\nonumber\\
&=& n^2n_1^2(n-n_1)^2 + \mid n_1(n-n_1)\mid(8n^2\mid n_1(n-n_1)\mid-6)\nonumber
\end{eqnarray*}

Using (\ref{n2}) we obtain that:
$$
R^2(n,n_1) \gtrsim  n^2n_1^2(n-n_1)^2
$$
By the same argument employed in \cite{Kenig1}, we state the following elemental
estimates without proof.
\begin{lemma}\label{log}
For any $\epsilon > 0$, $\alpha \in \mathbb{R}$ and $0 < \rho < 1$, we have:\\
$$\displaystyle{\int_{\mathbb{R}}}\frac {d\beta}{(1 +|\beta|)( 1 + |\alpha-\beta|)} \lesssim  \frac{\log(2+|\alpha|)}{(1+|\alpha|)}.$$
$$\displaystyle{\int_{\mathbb{R}}}\frac {d\beta}{(1 +|\beta|)^{\rho}( 1 + |\alpha-\beta|)} \lesssim  \frac{1+\log(1+|\alpha|)}{(1+|\alpha|)^{\rho}}.$$
$$\displaystyle{\int_{\mathbb{R}}}\frac {d\beta}{(1 +|\beta|)^{1+\epsilon}( 1 + |\alpha-\beta|)^{1+\epsilon}} \lesssim  \frac{1}{(1+|\alpha|)^{1+\epsilon}}.$$
\end{lemma}
\begin{lemma}\label{logc}
There exists $c > 0$ such that for any $\rho > \frac{2}{3}$ and any $\tau$, $\tau_1 \in \mathbb{R}$, the following is true :
$$
\displaystyle{\sum_{n_1\neq0}}\frac{\log( 2 + |\tau + m(n_1) + m(n-n_1)|)}{( 1 + |\tau + m(n_1) + m(n-n_1)|)} \leq C.
$$
$$
\displaystyle{\sum_{n\neq0}}\frac{\log( 2 + |\tau_1 + m(n_1) - m(n-n_1)|)}{( 1 + |\tau_1 + m(n_1) - m(n-n_1)|)} \leq C.
$$
$$
\displaystyle{\sum_{n\neq0}}\frac{\log( 1 + |\tau_1 + m(n_1) - m(n-n_1)|)}{( 1 + |\tau_1 + m(n_1) - m(n-n_1)|)^{\rho}} \leq C.
$$
\end{lemma}

\begin{proposition}\label{estimationbilinieaire}
Let $s \geq -\frac{1}{2}$, then for all $f$, $g$ with compact support in time included in the subset
 $\{(t,x), t \in [-T,T]\}$, there exists $\theta > 0 $ such that:
$$\lVert \partial _x(fg)\rVert_{X^{s,-\frac{1}{2}}} \lesssim T^{\theta} \lVert f \rVert_{X^{s,\frac{1}{2}}}\lVert g \rVert_{X^{s,\frac{1}{2}}}.
$$
\end{proposition}
\begin{remark}\label{lemmenonexistence}
This proposition is false for $s<-\frac{1}{2}$. We can exhibit
   a counterexample to the bilinear estimate in the Prop (\ref{estimationbilinieaire}) inspired by the similar argument in \cite{Kenig1}.


\end{remark}
 We now use the lower bound of the resonance function to recover the derivative on the non-linear term $u\partial_xu$.
\begin{lemma}\label{lower}
  Let 
$$
F_s = \frac{\mid n \mid^{2s+2}\mid n_1(n-n_1)\mid^{-2s}}{\sigma(\tau,\tau_1,n,n_1)}
$$
and 
$$
F_{s,r}=\frac{\mid n \mid^{2s+2}\mid n_1(n-n_1)\mid^{-2s}}{\sigma^{2(1-r)}(\tau,\tau_1,n,n_1)}
$$
where $\sigma(\tau,\tau_1,n,n_1) = \text{max}\{\mid \tau + m(n)\mid, \mid \tau_1 + m(n_1)\mid, \mid \tau - \tau_1 +  m(n-n_1)\mid\}$.
Then, for $ s \geq - \frac{1}{2}$, $0 < r < \frac{1}{4}$, we have 
$$ 
F_s \lesssim 1.
$$
and 
$$
F_{s,r} \lesssim \frac{1}{\mid n\mid ^{2-4r}}.
$$
\end{lemma}
\textbf{Proof}: This follows from  Lemma \ref{lemma1}.\\
\\
According to \cite{Ge95} we have the following Lemma: 
\begin{lemma}\label{sortirT1}
For any $ u \in X^{s,\frac{1}{2}}$ supported in $[-T,T]$ and for any $0< b < \frac{1}{2}$, 
 it holds:
\begin{equation}
 ||u||_{X^{s,b}} \lesssim  T^{(\frac{1}{2}-b)-}||u||_{X^{s,1/2-}}\lesssim  T^{(\frac{1}{2}-b)-}||u||_{X^{s,1/2}}.
\end{equation}
\end{lemma}

\textbf{Proof {of Proposition \ref{estimationbilinieaire}}} : Let$$P^b_f(n,\tau)= |n|^s < \tau + m(n) >^{b} |\hat{f}(n,\tau)|,$$
then we have 
$$
\|f\|_{X^{b,s}} = (\displaystyle{\sum_{n} \int_{\mathbb{R}}( P^b_f(n,\tau))^2d \tau)^{\frac{1}{2}} = \| P^b_f(n,\tau)\|_{l^2_nL^2_\tau}},
$$
and

\begin{equation}\label{B(fg)}
B(f,g)(n,\tau) = n^{s+1} < \tau + m(n) >^{-\frac{1}{2}}\displaystyle{\sum_{n_1\neq0,n_1\neq n}}\int_{\mathbb{R}}
\frac{(n_1(n-n_1))^{-s}P^{\frac{1}{2}-\gamma}_f(n_1,\tau_1)P^{\frac{1}{2}}_g(n-n_1,\tau-\tau_1)d\tau_1}{< \tau_1 + m(n_1) >^{\frac{1}{2}-\gamma} <\tau- \tau_1 + m(n-n_1) >^{\frac{1}{2}}}
\end{equation}
Denote $$F(n,\tau,n_1,\tau_1) = \frac{\mid n \mid^{s+1}\mid n_1(n-n_1)\mid^{-s}}{< \tau + m(n) >^{\frac{1}{2}}< \tau_1 + m(n_1) >^{\frac{1}{2}-\gamma} 
<\tau- \tau_1 + m(n-n_1) >^{\frac{1}{2}}}.
$$
Letting 
$E = \{(n,\tau,n_1,\tau_1): \mid\tau- \tau_1 + m(n-n_1)\mid \leq  \mid \tau_1 + m(n_1) \mid \},$
then by symmetry, (\ref{B(fg)}) is reduced to estimate
\begin{equation}\label{B(fg)1}
(\sum_{n\neq 0}\int_{\mathbb{R}}(\sum_{n_1\neq n,n_1\neq 0}\int_{\mathbb{R}} (1_E F)(n,\tau,n_1,\tau_1)P^{\frac{1}{2}-\gamma}_f(n-n_1,\tau-\tau_1)P^{\frac{1}{2}}_g(n_1,\tau_1)d\tau_1)^2d\tau)^{\frac{1}{2}}.	 
\end{equation}
We separate the two cases.\\
Case I:$\mid \tau_1 + m(n_1) \mid \leq \mid \tau + m(n) \mid$\\
In this case, the set E is replaced by
$$
E_I = \{(n,\tau,n_1,\tau_1): \mid \tau- \tau_1 + m(n-n_1)\mid \leq  \mid \tau_1 + m(n_1)  \mid \leq \mid \tau + m(n) \mid \},
$$
then by Cauchy-Schwarz inequality (\ref{B(fg)1}) is controled by
\begin{align}\label{B(fg)2}
\bigg\| \bigg(\sum_{n_1\neq n,n_1\neq 0}\int_{\mathbb{R}}(1_{E_I} &F)^2(n,\tau,n_1,\tau_1) d\tau_1\bigg)^{\frac{1}{2}}\nonumber \\
&\times \bigg(\sum_{n_1\neq n,n_1\neq 0}\int_{\mathbb{R}}
(P^{\frac{1}{2}-\gamma}_f)^2(n-n_1,\tau-\tau_1)(P^{\frac{1}{2}}_g)^2(n_1,\tau_1)d\tau_1\bigg)^{\frac{1}{2}}\bigg\|_{l^2_nL^2_\tau}.
\end{align}

Remark that $$
F^2 \approx F_s \frac{1}{< \tau_1 + m(n_1) >^{1-2\gamma} <\tau- \tau_1 + m(n-n_1) >},
$$
with $F_s = \frac{\mid n \mid^{2s+2}\mid n_1(n-n_1)\mid^{-2s}}{\sigma(\tau,\tau_1,n,n_1)}$, 
then by Lemma \ref{lower}, for $s \geq -\frac{1}{2}$, $(n,\tau,n_1,\tau_1) \in E_I$, we have
$$
\displaystyle{\sup_{n,\tau}\sum_{n_1}\int_{\mathbb{R}}(1_{E_I} F)^2(n,\tau,n_1,\tau_1)d\tau_1 
\lesssim \sup_{n,\tau}\sum_{n_1}}\int_{\mathbb{R}} \frac{d\tau_1}{< \tau_1 + m(n_1) >^{1-2\gamma} <\tau- \tau_1 + m(n-n_1) >}
$$
we can easily see that
$$
(\ref{B(fg)2})  \leq \sup_{n,\tau}\sum_{n_1}\int_{\mathbb{R}} \frac{d\tau_1}{< \tau_1 + m(n_1) >^{1-2\gamma} <\tau- \tau_1 + m(n-n_1) >}
\lVert P^{\frac{1}{2}-\gamma}_f(n,\tau) \lVert_{l^2_nL^2_\tau}\lVert P^{\frac{1}{2}}_g(n,\tau) \lVert_{l^2_nL^2_\tau} 
$$
then by Lemma \ref{log}, \ref{logc}( take $\alpha = \tau + m(n_1) +m(n-n_1)$ and $\beta = \tau_1 + m(n_1)$) and \ref{sortirT1} we obtain that there exist
$\theta >0$ such that:
$$
(\ref{B(fg)1}) \lesssim \lVert f \rVert_{X^{s,\frac{1}{2}-\gamma}}\lVert g \rVert_{X^{s,\frac{1}{2}}}\lesssim T^{\theta}
\lVert f\rVert_{X^{s,\frac{1}{2}}}\lVert g \rVert_{X^{s,\frac{1}{2}}}.
$$
Case II:$  \mid \tau + m(n) \mid \leq \mid \tau_1 + m(n_1) \mid $
Here the set E becomes:
$$
E_{II} = \{(n,\tau,n_1,\tau_1): \mid\tau- \tau_1 + m(n-n_1)\mid \leq  \mid \tau_1 + m(n_1) \mid, \mid \tau + m(n) \mid < \mid \tau_1 + m(n_1) \mid\}.
$$
Then we will estimate
\begin{equation}\label{B(fg)4}
\lVert\sum_{n_1}\int_{\mathbb{R}}(1_{E_{II}} F)(n,\tau,n_1,\tau_1)P^{\frac{1}{2}-\gamma}_f(n-n_1,\tau-\tau_1)P^{\frac{1}{2}}_g(n_1,\tau_1)d\tau_1\parallel_{l^2_nL^2_\tau}
\end{equation}
By duality, (\ref{B(fg)4}) equals to 
\begin{equation}\label{B(fg)5}
 \displaystyle{\sup_{\lVert w \rVert_{l^2_nL^2_\tau} = 1}\sum_{n,n_1}\int_{\mathbb{R}^{2}}
 w(n,\tau)(1_{E_{II}} F)(n,\tau,n_1,\tau_1)P^{\frac{1}{2}-\gamma}_f(n-n_1,\tau-\tau_1)P^{\frac{1}{2}}_g(n_1,\tau_1)d\tau_1d\tau}.
\end{equation}
By Fubini's Theorem and Cauchy-Schwarz inequality, we could control (\ref{B(fg)5}) by 
\begin{align}\label{B(fg)6}
\sup_{\lVert w \rVert_{l^2_nL^2_\tau} = 1}\bigg(\sum_{n_1}&\int_{\mathbb{R}}\big[\sum_n \int_{\mathbb{R}}(1_{E_{II}} F)^2(n,\tau,n_1,\tau_1)d\tau\big]
\times\\
&\big[\sum_{n}\int_{\mathbb{R}} w^2 (P^{\frac{1}{2}-\gamma}_f)^2(n-n_1,\tau-\tau_1)d\tau\big]d\tau_1\bigg)^{\frac{1}{2}}\lVert g \rVert_{X^{s,\frac{1}{2}}}.
\nonumber
\end{align}
Similary to the previous case, we can show that:
$$
\displaystyle{\sup_{n_1,\tau_1} \sum_n \int_{\mathbb{R}}(1_{E_{II}} F)^2(n,\tau,n_1,\tau_1)d\tau \lesssim 1}.
$$
Finaly we obtain that
$$
(\ref{B(fg)6})\lesssim
 \lVert f \rVert_{X^{s,\frac{1}{2}-\gamma}} \lVert g \rVert_{X^{s,\frac{1}{2}}} \lesssim T^{\theta}
\lVert f \rVert_{X^{s,\frac{1}{2}}}\lVert g \rVert_{X^{s,\frac{1}{2}}}.
$$


Now we have  the following proposition:
\begin{proposition}\label{l1l2}
Let $s \geq -\frac{1}{2}$ then for all  $f$, $g$ with compact support in time included in the subset
 $\{(t,x), t \in [-T,T]\}$, there exists $\theta> 0 $ such that:
 \begin{equation}\label{l1l2bis}
\bigg(\sum_{n \in \mathbb{\dot{Z}}}\mid n\mid^{2s}\bigg[\int_{\mathbb{R}}
\frac{\mid n \hat{f}*\hat{g}(n,\tau)\mid}{ <\tau + m(n) >}d\tau\bigg]^2 \bigg)^{\frac{1}{2}} \lesssim T^{\theta}
\lVert f \rVert_{X^{s,\frac{1}{2}}}\lVert g \rVert_{X^{s,\frac{1}{2}}}.
\end{equation}
\end{proposition}
\textbf{Proof:} As in the proof of Prop \ref{estimationbilinieaire}, we consider (\ref{l1l2bis}) in the same two cases. It could be written as:

\begin{equation}\label{l1l2one}
\bigg\| \int_{\mathbb{R}}\sum_{n_1}\int_{\mathbb{R}} (1_E F)(.,\tau,n_1,\tau_1)P_f^{\frac{1}{2}-\gamma}(.-n_1,\tau-\tau_1)P_g^{\frac{1}{2}}(n_1,\tau_1)
d\tau_1d\tau \bigg\|_{l^2_n} \lesssim T^{\theta }\lVert f\rVert_{X^{s,\frac{1}{2}}}\| g\|_{X^{s,\frac{1}{2}}},	 
\end{equation}
where
$$
F(n,\tau,n_1,\tau_1) = \frac{\mid n \mid^{s+1}\mid n_1(n-n_1)\mid^{-s}}{< \tau + m(n) >^{\frac{1}{2}}< \tau_1 + m(n_1) >^{\frac{1}{2}-\gamma} 
<\tau- \tau_1 + m(n-n_1) >^{\frac{1}{2}}}.
$$
1)Case I: $\mid \tau_1 + m(n_1) \mid \leq \mid \tau + m(n)\mid$.
As before, the set $E$ is replaced by 
$$
E_{I}=\{\mid \tau- \tau_1 + m(n-n_1) \mid \leq \mid  \tau_1 + m(n_1) \mid \leq \mid \tau + m(n)\mid\}.
$$
By duality 
, we suffer to estimate 
$$
\displaystyle{\sup_{\lVert w \rVert_{l^2_n} = 1}\sum_{n,n_1}\int_{\mathbb{R}^{2}}
 w(n)(1_{E_{I}} F)(n,\tau,n_1,\tau_1)P^{\frac{1}{2}-\gamma}_f(n-n_1,\tau-\tau_1)P^{\frac{1}{2}}_g(n_1,\tau_1)d\tau_1d\tau}.
$$
Now by Cauchy-Schwarz, we could control it by
\begin{align}
\sup_{\lVert w \rVert_{l^2_n} = 1}\big(\sum_{n_1}\int_{\mathbb{R}}\big[\sum_n &\int_{\mathbb{R}}(1_{E_{I}} F)^2(n,\tau,n_1,\tau_1)d\tau\big]\times
\nonumber\\
&\big[\sum_{n}\int_{\mathbb{R}} w^2 (P^{\frac{1}{2}-\gamma}_f)^2(n-n_1,\tau-\tau_1)d\tau\big]d\tau_1\big)^{\frac{1}{2}}\lVert g \rVert_{X^{s,\frac{1}{2}}},
\nonumber
\end{align}
then it is sufficient to show that, for $s \geqslant -\frac{1}{2}$
$$
D=\displaystyle{\sup_{n_1} \sum_n \int_{\mathbb{R}}\int_{\mathbb{R}}(1_{E_{I}} F)^2(n,\tau,n_1,\tau_1)d\tau d\tau_1\lesssim 1}.
$$
For some $0 < r < \frac{1}{4}$, D can be rewriten as:
$$
D=\displaystyle{\sup_{n_1} \sum_{n} \int_{\mathbb{R}}\int_{\mathbb{R}}\frac{1}{< \tau + m(n) >^{2r}},
(1_{E_{I}} F_r)^2(n,\tau,n_1,\tau_1)d\tau d\tau_1}
$$
where 
$$
F_r^2 = \frac{\mid n \mid^{2s+2}\mid n_1(n-n_1)\mid^{-2s}}{< \tau + m(n) >^{2(1-r)}}\frac{1}{< \tau_1 + m(n_1) >^{1-2\gamma} <\tau- \tau_1 + m(n-n_1) >}.
$$
Remark that
$$F_r^2= F_{s,r}\frac{1}{< \tau_1 + m(n_1) >^{1-2\gamma} <\tau- \tau_1 + m(n-n_1)>},$$
then by Lemma \ref{lower}, D could be controlled by
$$
D \lesssim \displaystyle{\sup_{n_1}\sum_{n}\int_{\mathbb{R}}\int_{\mathbb{R}}\frac{1}{\mid n \mid ^{2-4r}}\frac{d\tau_1d\tau}{<\tau_1 + m (n_1)>^{1-2\gamma+r} 
< \tau- \tau_1 + m(n-n_1) >^{1+r}}}
$$ 
 by Lemma \ref{log} we have:
$$
\int_{\mathbb{R}}\frac{d\tau_1}{<\tau_1 + m (n_1)>^{1-2\gamma+r} 
< \tau- \tau_1 + m(n-n_1) >^{1+r}}\lesssim  \frac{1}{(1 + \mid \tau + m(n_1) +m(n-n_1)\mid)^{1-2\gamma+r}}.
$$
Hence
$$
D \lesssim \displaystyle{\sup_{n_1}\sum_{n}\frac{1}{\mid n \mid ^{2-4r}}\int_{\mathbb{R}} \frac{d\tau}{(1 + \mid \tau + m(n_1) +m(n-n_1)\mid)^{1-2\gamma+r}}}.
$$
Therefore, if $ r < \frac{1}{4}$, we have $D \lesssim \sum_n\frac{1}{\mid n \mid ^{2-4r}} < +\infty$.\\
2)Case II, $  \mid \tau + m(n) \mid \leq \mid \tau_1 + m(n_1) \mid $.
Now we replace $E$ with
$$
E_{II} = \{(n,\tau,n_1,\tau_1): \mid\tau- \tau_1 + m(n-n_1)\mid \leq  \mid \tau_1 + m(n_1) \mid, \mid \tau + m(n) \mid < \mid \tau_1 + m(n_1) \mid\}.
$$
We write 
$$
1+\mid \tau + m(n)\mid = (1+\mid\tau+m(n)\mid)^r(1+\mid\tau + m(n)\mid)^{1-r},
$$
where $\frac{1}{2} < r < 1$. As in case I, $F_r$ denotes
$$
\frac{\mid n \mid^{s+1}\mid n_1(n-n_1)\mid^{-s}}{< \tau + m(n) >^{(1-r)}}\frac{1}{<\tau_1 + m(n_1)>^{\frac{1}{2}-\gamma}
<\tau- \tau_1 + m(n-n_1) >^{\frac{1}{2}}}.
$$
It suffices to estimate 
\begin{equation}\label{E2Fr}
\Bigg(\sum_n\bigg(\int_{\mathbb{R}}\sum_{n_1}\int_{\mathbb{R}} \frac{1}{< \tau + m(n)>^r}(1_{E_{II}}F_r)(n,\tau,n_1,\tau_1)
P^{\frac{1}{2}-\gamma}_f(n-n_1,\tau-\tau_1)
P^{\frac{1}{2}}_g(n_1,\tau_1)d\tau_1 d\tau \bigg)^2\bigg)^{\frac{1}{2}}.	 
\end{equation}
Applying the Cauchy-Schwarz inequality in $\tau$ we see that (\ref{E2Fr}) is bounded by 
\begin{align}
\displaystyle{\Bigg[\sum_{n}\bigg(\int_{\mathbb{R}}}& \frac{d\tau}{< \tau + m(n)>^{2r}}\bigg)\nonumber\\
&\times\bigg(\int_{\mathbb{R}}\bigg(\sum_{n_1}\int_{\mathbb{R}}
(1_{E_{II}}F_r)(n,\tau,n_1,\tau_1)
P^{\frac{1}{2}-\gamma}_f(n-n_1,\tau-\tau_1)
P^{\frac{1}{2}}_g(n_1,\tau_1)d\tau_1 \bigg)^2 d\tau\bigg)\Bigg]^{\frac{1}{2}}.\nonumber
\end{align}
Since $2r > 1$, then (\ref{E2Fr}) is dominated by 

$$\displaystyle{\Bigg[\sum_{n}\int_{\mathbb{R}}
\bigg(\sum_{n_1}\int_{\mathbb{R}}
(1_{E_{II}}F_r)(n,\tau,n_1,\tau_1)
P^{\frac{1}{2}-\gamma}_f(n-n_1,\tau-\tau_1)
P^{\frac{1}{2}}_g(n_1,\tau_1)d\tau_1 \bigg)^2 d\tau\bigg)\Bigg]^{\frac{1}{2}}},$$
then as the case II in the proof of Prop \ref{estimationbilinieaire}, we obtain the estimate, and this end the proof.


Now we return to the proof of \textbf{Theorem \ref{thexistence}}: Let $L$ defined by

\begin{equation}\label{FormDuhTrOst}
L(u) = \psi(t)[S(t)\phi - \int_0^t S(t-t')\partial_x(\psi_T^2u^2(t'))dt'],
\end{equation}
where $t\in \mathbb R$, $\psi$ indicates a time cutoff function~:
\begin{equation}\label{Cutoff}
\psi \in C_0^{\infty}(\mathbb R),\quad \text{sup }\psi\subset [-2,2], \quad \psi=1\text{ on }[-1,1],
\end{equation}
$ \psi_T(.)=\psi(./T),$\\
we will apply a fixed  point argument to (\ref{FormDuhTrOst}), using the following estimates:
\begin{proposition}\label{contraction}
There exists a constant $C = C(\phi)$ such that:
$$
||L(u)||_{Y^{s}}\leq C||\phi||_{H^{s}}+
CT^{\gamma}||u||^2_{Y^{s}}
$$
and 
$$
||L(u)-L(v)||_{Y^{s}}\leq CT^{\gamma}||u-v||_{Y^{s}}||u+v||_{Y^{s}}.
$$
 
\end{proposition}
Proof: It follows from Propositions \ref{estimationbilinieaire}, \ref{l1l2} and classical linear estimates (see \cite{Colliander}).\\

Note that if we take  $T=(4C^2||\phi||_{H^{s}})^{-1/\gamma}$ we deduce
from Prop \ref{contraction} that $L$ is strictly contractive in  the ball $B(0,\frac{1}{8C^2})$ in $Y^{s}$. This proves the existence of 
a unique solution $u$  to (\ref{FormDuhTrOst}) in $Y^{s}$.
\subsection{Global existence in $L^2$}
Its easy to see that the $L^2$-norm is conserved ( $\|u(t)\|_{L^2} = \|u_0\|_{L^2}$). Hence, if we take an initial data $u_0$ in $L_0^2$, 
 the solution $ u $ such that $u(0)=u_0$ can be extended for all positive times and the existence is global in $L_0^2(\mathbb{T})$.

\end{document}